\documentclass[11pt]{elsarticle}
\usepackage{amsmath,amssymb,amsthm, xcolor,url}

\usepackage{makeidx}
\usepackage{graphicx, natbib}
\usepackage[a4paper]{geometry}
\usepackage[displaymath]{lineno}
\usepackage[singlespacing]{setspace}%
\usepackage[colorlinks=true,breaklinks=true,bookmarks=true,urlcolor=blue,
     citecolor=blue,linkcolor=blue,bookmarksopen=false,draft=false]{hyperref}
     \definecolor{DarkOrchid} {cmyk}{0.40,0.80,0.20,0}
\def\tcb{\textcolor{blue}}

  \biboptions{sort&compress}
\providecommand{\U}[1]{\protect\rule{.1in}{.1in}}

\def\dom {\mathop{\rm Dom\,  }}
\newtheorem{theorem}{Theorem}

\newtheorem{corollary}[theorem]{Corollary}

\newtheorem{definition}[theorem]{Definition}
\newtheorem{remark}[theorem]{Remark}

\def\R{\Bbb R}
\def\N{\Bbb N}

\def\dom {\mathop{\rm dom\,  }}

\def\gph{ \mathop{\rm gph\,  }}

\begin{document}

\title{ Ekeland's  inverse function theorem in graded Fr\'{e}chet spaces  revisited for multifunctions\tnoteref{t1}  }
\tnotetext[t1]{Dedicated to the memory of Jonathan M. Borwein.}
\author{Huynh Van Ngai\fnref{fn1}}
\ead{nghiakhiem@yahoo.com}
\address{Department of Mathematics, University of
Quynhon}

\author{Michel  Th\'era\corref{cor2}\fnref{fn2}}
\ead{michel.thera@unilim.fr}
\address{Universit\'{e} de Limoges, Laboratoire XLIM, UMR-CNRS 6172, France   and Centre for Informatics and Applied Optimization, Federation University, Australia }
\cortext[cor2]{Principal corresponding author} 
\fntext[fn1]{Research supported by NAFOSTED under Grant 101.01-2016.27 and by Vietnam Institute for Avanced Study in Mathematics (VIASM).}    
\fntext[fn2]{Research supported by  the Gaspard Monge Program for Optimization and  by
 Australian Research Council, project  DP160100854.}
\begin{abstract}
 In this paper, we present some {inverse} function theorems {and implicit function theorem} for set-valued mappings between Fr\'echet spaces. The proof relies on Lebesgue's Dominated Convergence Theorem and on Ekeland's variational principle. An application to the existence of solutions of differential equations in Fr\'{e}chet spaces with non-smooth data is given.

\end{abstract}


\maketitle

\section{Introduction} 
 The {inverse} function theorem is {one of the central components } of the classical {and} the modern variational analysis and  {an essential} device to solving nonlinear equations.  The {inverse} function theorem or its variants known as the {implicit} function theorem  or the rank theorem have  been established originally in Euclidean spaces  and then extended to the Banach space setting. Outside this setting, for instance in  Fr\'echet spaces, it is known that the inverse function theorem generally fails 
  (see  
Lojasiewicz Jr \& 
Zehnder \cite{LZ}).
This is the reason why another form of   {inverse} function theorem,  nowadays called the Nash-Moser theorem is used  as a powerful tool 
to prove local existence for non-linear partial differential equations in spaces of smooth functions.
Some inverse theorems  of  Nash-Moser type have also  been proved for {mappings} between Fr\'echet spaces,    { that are  supposed to be  \textit{tame}, an additional  property  guaranteeing  that the semi-norms satisfy some interpolation
properties, see e.g. \cite{Hamilton, Gerard})}  or {that allow the  use of  smoothing operators}  as introduced by Nash (see e.g. \cite{Hamilton, Hor, LZ, Moser, Sergeraert, Zenhder, Schwartz}). {To overcome the loss of derivatives, these additional properties in Fr\'echet spaces allow Newton's method on which   the Nash-Moser type inverse function theorems   are based to converge.}  
Recently, Ekeland \cite{E-IHP}  (see also Ekeland \& S\'er\'e \cite{E-S}) produced a new result within a class  {of spaces} much larger than the
one used in the Nash-Moser literature.
\vskip 2mm
Nowadays, modeling has evolved beyond equations and we know the importance and the efficacy of studying 
set-valued solution mappings which
assign to each instance of the parameter element in a model all the corresponding
solutions, if any.  As it is mentioned in the book by Dontchev \& Rockafeller \cite{DR},  ``\textit{the central question is whether a {solution mapping} can be localized
graphically in order to achieve single-valuedness and in that sense produce a
function, the desired implicit function}". To be more explicit,
many applied problems can be modeled as differential inclusions or more generally as generalized equations, that is, inclusions governed by a set-valued mapping. For these problems which are the analogous of nonlinear equations, there is a need to use implicit multifunction theorems.  During the last years a wide literature has emerged  related to implicit multifunction theorems (see e.g., \cite{durea2012, NT, HNT, Huu, Nghia, Pang, AB, Ledyaev, DR} and the references therein). However, to our knowledge, they have been established in the framework of Banach spaces and nothing exists for Fr\'echet spaces. Therefore, motivated by the recent work by  Ekeland, and Ekeland \& S\'er\'e,  
it is our  aim in this paper to investigate the possibility  to obtain, in the context of graded Fr\'echet spaces,  an implicit multifunction theorem for set-valued mappings.
\vskip 2mm
{The structure of the paper is as follows. Section 1 is devoted to recalling the notions of Fr\'echet, graded Fr\'echet, standard Fr\'echet spaces, and contingent derivative of multifunctions, concepts essential for the framework  and assumptions on which our results are based.  In Section 2, we present an inverse multifunction theorem for set-valued mappings between Fr\'echet spaces which is the main result of the paper. This result allows us to obtain a general version of the Ekeland inverse function for G\^{a}teaux differentiable mappings and to establish an implicit multifunction theorem for parametrized set-valued mappings. In the final section, we present an application to the existence of solutions for differential equations in   Fr\'{e}chet spaces. } 
\section{Preliminaries}
We begin this section  with recalling  {briefly} some notions on Fr\'{e}chet spaces, i.e. on  locally convex spaces which are Hausdorff, complete and whose topology is induced by a countable family of semi-norms $(\Vert\cdot\Vert_{k})_{k\in \N}$ {with the property:
$$x\in X, \quad \;\|x\|_k=0\;\text{for all}\; k\in\N\quad \Rightarrow\; x=0.$$}  { This class of spaces contains evidently Banach spaces, as well as many other locally convex spaces used in various  areas of real or complex analysis.}
{Also notice that given a Fr\'echet space $F$, we may produce the Fr\'echet space  $C^0([a,b], F)$ of continuous paths in $F$ equipped with the semi-norms defined by $\Vert f\Vert_k:=\sup_{t\in [a,b]} \Vert f(t)\Vert_k$ and widely used in analysis.}

A  \textit{graded Fr\'echet } space $X$ is a Fr\'{e}chet space,    whose topology is generated by a fixed 
sequence of semi-norms $(\|\cdot\|_k)_{k\in\N}$, increasing in strength, so that, 
$$\|x\|_k\le\|x\|_{k+1}\quad\forall x\in X,\;\forall k\in\N.$$
{ This class contains}  the space $C^\infty ([a,b])$ of infinitely differentiable {real-valued} functions on the interval $[a,b]$ with the grading 
$$\Vert f\Vert_n= \sup_{k\leq n}\sup_{x\in [a,b]}\vert D^k(f(x))\vert,$$
as well as the space  $C^\infty(\overline{\Omega}, \R^d)$, 
where $\overline{\Omega}\subset \R^n$ is compact, with a smooth boundary and is  the
closure of its interior are  graded Fr\'{e}chet spaces. { Note that  every Fr\'echet space can be considered as a  graded Fr\'echet  space by replacing the initial  family of semi-norms by the semi-norms $\Vert\cdot \Vert_{n}:= \displaystyle\sum_{k=1}^{k=n}\Vert\cdot\Vert_k.$ However, as the Nash-Moser inverse function theorem highlights, 
 the grading  plays  an essential role in its statement and in its proof.}

{ Well known are the facts:
\begin{itemize}
\item    the Cartesian product of two graded Fr\'echet spaces is a graded Fr\'echet space with the grading $\Vert (x,y)\Vert_n= \Vert x\Vert_n+\Vert y\Vert_n;$ 
\item  a closed subspace of a graded Fr\'echet space is a graded Fr\'echet space;
\item   a sequence $(x_n)_{n\in \N}$ of  elements in  a Fr\'echet space $X$ converges to $x\in X$,  if and only if,  $\|x_n-x\|_k\to 0$ for all $k\ge 0;$  
\item  due to a classical result that (see e.g. \cite{E-IHP}), a graded Fr\'{e}chet space is a complete metric space with the metric:
\begin{equation}\label{Dist}
d(x,y)=\sum_{k=1}^\infty\mu_k\min\{r,\|x-y\|_k\},
\end{equation}
where, $(\mu_k)_{k\in \N}$ is any sequence of non-negative numbers with unbounded support: $\sup\{k\in\N:\; \mu_k\not=0\}=+\infty$ and $r>0.$ 
\end{itemize}
}
\vskip 0.2cm
\begin{definition}\label{Def-Standard} (\cite{E-IHP}, Definition 5)  A graded Fr\'{e}chet space is said to be standard if for every $x\in X,$ we can find a constant {$c:=c(x)$} a sequence $(x_n)\subseteq X$ converging to $x$ and a sequence of non-negative  numbers {$(c_n)$}  such that
\begin{equation}\label{Stand}
\|x_n\|_k\le {c}\|x\|_k\quad \mbox{and}\quad \|x_n\|_k\le {c_n}^k\quad\forall k,n\in\N.
\end{equation}
\end{definition}
The space  $C^\infty(\overline{\Omega}, \R^d)$ is a standard graded Fr\'echet space, see \cite{E-IHP}.
\vskip 0.2cm
Let $X$, $Y$ be graded Fr\'{e}chet spaces. Consider a multifunction (set-valued mapping)  $F:X\rightrightarrows Y$ between $X$ and $Y$, that is a function between $X$ and the subsets (possibly empty) of $Y$.  We denote  by $\gph F$ and $\dom F$ and $F^{-1}: Y\rightrightarrows X,$ the {\em graph,} {\em domain} and {\em inverse } of $F,$ respectively:  $$\gph F=\{(x,y)\in X\times Y:\; y\in F(x)\};\;\; \dom F=\{x\in X:\; F(x)\not=\emptyset\}$$
and
$$F^{-1}(y):=\{x\in X:\;\; y\in F(x)\},\quad y\in Y.$$
We say that $F$ is a {\em closed} multifunction if $\gph F$ is a closed subset of $X\times Y.$ In what follows, we will use the notion of  \textit{ contingent derivative} of multifunctions. The contingent derivative of the multifunction $F$ at a point $(\bar x,\bar y)\in \gph F$ is the  multifunction  $DF(\bar x,\bar y): X\rightrightarrows Y$, defined  for $u\in X$ by 
$$\{v\in Y:\;\exists (t_n)_{n\in\N}\downarrow  0^+, \; \exists (u_n, v_n)_{n\in\N}\to (u, v)\;\mbox{with}\; (\bar x+t_nu_n,\bar y+t_nv_n)\in\gph F,\;\forall n\}.$$
In other words,
$DF(\bar x,\bar y)(u)=\{v\in Y:\;\; (u,v)\in T_{\gph F}(\bar x,\bar y)\},$ where, $T_{\gph F}(\bar x,\bar y)$ stands for the \textit{{contingent cone}} to $\gph F$ at $(\bar x,\bar y).$ For more details, the reader is referred to  the book by Aubin \& Frankowska \cite{AF}.

 When $F: X\rightarrow Y$ is a single-valued mapping, we use the notation  $DF(\bar x)$  for $DF(\bar x,\bar y)$. Note that if $F:X\to Y$ is G\^{a}taux differentiable at $\bar x\in X$, then $DF(\bar x)$ coincides with  the G\^ateaux derivative of $F$ at $\bar x.$

\section{{Inverse and} Implicit multifunction theorems}
Throughout, we consider a closed multifunction $F:X\rightrightarrows Y$ between graded Fr\'{e}chet spaces $X$ and $Y$ {induced respectively by countable families of semi-norms for which we use the same notation $(\|\cdot\|_k)_{k\in\N}.$ For each $k\in\N$, $x\in X$ and a subset $S\in X$, denote by $d_k(x,S):=\inf_{z\in X}\|z-x\|_k$, which is referred to {as} the {\it semi-distance} from the point $x$ to the set $S$  with respect to the semi-norm $\|\cdot\|_k$ in $X.$} For given $r\in (0,+\infty],$  $k_0\in\N,$  we denote  respectively by
$B_X(\bar x,k_0,r)=\{x\in X:\;\; \|x-\bar x\|_{k_0}<r\} $ and $B_X[\bar x,k_0,r]=\{x\in X:\;\; \|x-\bar x\|_{k_0}\le r\}$, {which are also referred to} {as}  the open and closed balls in $X$ centered at $\bar x\in X$,  with radius $r$ with respect to the semi-norm $ \|\cdot\|_{k_{0}}$. {Note that since $\|\cdot\|_{k_0}$ is a semi-norm, for any $z\in X$ with $\|z-\bar x\|_{k_0}=0,$ $B_X(z,k_0,r)$, $B_X[z,k_0,r]$ coincide with $B_X(\bar x,k_0,r)$, $B_X[\bar x,k_0,r]$, respectively.} 
\vskip 0.2cm
\begin{theorem}\label{MS}
Let $F:X\rightrightarrows Y$ be a closed multifunction between graded Fr\'{e}chet spaces and let  $(\bar x,\bar y)\in\gph F$ {be given}. Assume furthermore that $Y$ is standard. Suppose also that there are integers $k_0, d_1,d_2,$ real numbers $r\in(0,+\infty]$,  $C\ge 0$  and non-decreasing sequences of non-negative  reals $(\nu_k)_{k\in \N},$ $(\nu^\prime_k)_{k\in \N}, (m_k)_{k\in \N},$   and $({a}_k)_{k\in \N}$ with $m_k\ge 1,$ ${a}_k\ge 1$ such that  the following conditions are satisfied:
\begin{itemize}
\item[(i)]  For all $(x,y)\in \gph F$  with $x\in B_X(\bar x,k_0,r),$ $y\in B_Y(\bar y,k_0+d_1+d_2,2r/{a}_{k_0+d_1}),$ for every { $(u,v)\in \gph DF(x,y),$} there exist {$c_2(u,v)>0$} and  sequences $t_n\downarrow 0,$ $ u_n\to u$ and $v_n\to v$ with $(x+t_nu_n,z+t_nv_n)\in\gph F$, such that for all $n\in\N,$ all $k\in\N,$
$$\|v_n\|_k\le {c_2(u,v)}(m_k\|u\|_{k+d_1}+\|x-\bar x\|_k/{a}_{k-d_2}+\|y-\bar y\|_k+\nu_k) $$
and 
$$\|u_n\|_k\le {c_2(u,v)}(m_k\|u\|_{k+d_1}+\|x-\bar x\|_k+\nu_k) ;$$
\item[(ii)]  For all $(x,y)\in\gph F$ with $x\in B_X(\bar x,k_0,r),$ $y\in B_Y(\bar y,k_0+d_1+d_2,2r/{a}_{k_0+d_1}),$  for  every $v\in Y,$  there exists $u\in DF(x,y)^{-1}(v)$ such that
$$\|u\|_k\le C\left(\frac{\|x-\bar x\|_{k-d_1-d_2}}{m_k{a}_{k-d_1-d_2}}+\nu^\prime_k\right)\|v\|_{d_1+d_2}+{a}_k\|v\|_{k+d_2},\;\; \forall k\in\N.$$
By convention we set $\|\cdot\|_{k}=\|\cdot\|_{0}$ and ${a}_k=1$ for $k<0.$ 
\end{itemize}
Let $(\beta_k)_{k\in \N}$ be a sequence of non-negative reals with unbounded support such that
\begin{equation}\label{nu-cond}
\begin{array}{ll}&\sum_{k=0}^\infty\beta_k\nu_k<+\infty \quad \sum_{k=0}^\infty\beta_{k}m_k\nu^\prime_{k+d_1}<+\infty\\
&\text{and}\quad \sum_{k=0}^\infty\beta_km_k{a}_{k+d_1}n^k<+\infty.
\end{array}
\end{equation}Then, for every $y\in Y$ with $C\gamma<1$, where
\begin{equation}\label{Cond-image-bis}
\gamma:=\frac{\sum_{k=0}^\infty\left(\beta_k\|y-\bar y\|_k+\beta_{k+d_1+d_2}\nu^\prime_{k+d_1}/{a}_{k+d_1}\right)}{\sum_{k=d_1+d_2}^\infty\beta_k},
\end{equation}
and
\begin{equation}\label{Cond-Image}
\sum_{k=0}^\infty\beta_k\|y-\bar y\|_k\left(1-\sqrt{C\gamma}\right)^{-2}<\frac{r\beta_{k_0+d_1+d_2}}{{a}_{k_0+d_1}},
\end{equation}
there exists $x\in B_X(\bar x,k_0,r)$ such that $y\in F(x).$
\end{theorem}
\vskip 0.2cm
\textit{Proof.} The proof   is based on the Ekeland variational principle 
\cite{Eke74}.
 However, the function and space  to which the Ekeland variational principle is applied are different from \cite{E-IHP}.
\vskip 0.1cm 
By translation if necessary, we can assume without loss of generality  that $\bar x=0$ and $\bar y=0.$ Let $(\alpha_k)_{k\in \N}$  be  the  sequence  defined by
$$\alpha_k=\frac{\beta_{k+d_1+d_2}}{{a}_{k+d_1}},\;\; k\in\N.$$
Consider the distances on $X$ and $Y$  defined respectively by
$$d(x_1,x_2):=\sum_{k=0}^\infty\alpha_k\min\{r, \|x_1-x_2\|_k\},\;\; x_1,x_2\in X,$$
$$d(y_1,y_2):=\sum_{k=0}^\infty\beta_k\min\{r, \|y_1-y_2\|_k\},\;\; y_1,y_2\in Y.$$
For $\varepsilon>0,$ we define the distance  $d_\varepsilon (\cdot,\cdot)$ on $X\times Y$ by
$$d_\varepsilon((x_1,y_1),(x_2,y_2)):=d(x_1,x_2)+\varepsilon d(y_1,y_2),\;\; (x_1,y_1),(x_2,y_2)\in X\times Y.$$
Equipped with  these distances,  the spaces $X,Y$ and therefore $X\times Y$ are complete metric spaces.  Let $y_0\in Y$ be such that  (\ref{Cond-image-bis}) and (\ref{Cond-Image}) are satisfied. {Setting $\eta=\sqrt{C\gamma}-C\gamma,$ where 
$\gamma$ is defined by (\ref{Cond-image-bis}), }
consider the extended-real-valued function $f:X\times Y\rightarrow \R\cup\{+\infty\}$ defined by
\begin{equation}\label{Func}
f(x,y)=\eta\sum_{k=0}^\infty\alpha_k\|x\|_{k+d_1}+\sum_{k=0}^\infty\beta_k\|y-y_0\|_k+\delta_{\gph F}(x,y), \;\; (x,y)\in X\times Y,
\end{equation} 
where $\delta_{\gph F}$ stands for the \textit{indicator function}  of $\gph F,$ that is,
$$\delta_{\gph F}(x,y)=\left\{\begin{array}{ll}0\quad&\mbox{if}\; (x,y)\in\gph F,\\
+\infty\quad&\mbox{otherwise}.
\end{array}\right.$$
\noindent \textit{\textbf{Claim 1. } The function $f$ is lower semicontinuous and  bounded from below on $X\times Y.$} 
\vskip 2mm
One has
\begin{equation}
0\le \inf_{(x,y)\in X\times Y}f(x,y)\le f(0,0)=\sum_{k=0}^\infty\beta_k\|y_0\|_k<+\infty.
\end{equation}
Take a sequence  $((x_n,y_n))_{n\in \N}$ converging to $(x,y)$ in $X\times Y.$ Then, $\|x_n-x\|_k\to 0$ and $\|y_n-y\|_k\to0$ for every $k\in \N.$ Two cases may happen:
\\
1.-  If $(x,y)\notin\gph F,$ then by the closedness of the graph of $F,$ $(x_n,y_n)\notin\gph F$ when $n$ is sufficiently large. Hence, 
$$\lim_{n}f(x_n,y_n)=+\infty=f(x,y).$$
2.-
 Suppose now that  $(x,y)\in\gph F.$  Thanks to the Fatou lemma we have,
$$\begin{array}{ll}
\displaystyle\liminf_{n\to +\infty}f(x_n,y_n)&\ge \eta\displaystyle\liminf_{n\to +\infty}\left(\sum_{k=0}^\infty\alpha_k\|x_n\|_k+\sum_{k=0}^\infty\beta_k\|y_n-y_0\|_k\right)\\
&\ge \eta\sum_{k=0}^\infty\alpha_k\displaystyle\lim_{n\to +\infty}\|x_n\|_k+\sum_{k=0}^\infty\displaystyle\lim_{n\to +\infty}\beta_k\|y_n-y_0\|_k\\
&=\eta\sum_{k=0}^\infty\alpha_k\|x\|_k+\sum_{k=0}^\infty\beta_k\|y-y_0\|_k=f(x,y),
\end{array}$$
establishing the 
claim.\qed
\vskip 0.2cm
As
$$\frac{C\gamma+\eta}{\eta(1-C\gamma-\eta)}=(1-\sqrt{C\gamma})^{-2},$$
in view of assumption (\ref{Cond-Image}),  take $\bar r>0$ and $\varepsilon>0$ such that
\begin{center}
\begin{equation}\label{Para1}
\bar r<r \quad \text{and} \quad \frac{C\gamma+\eta+\eta\varepsilon}{\eta(1-C\gamma-\eta)}\sum_{k=0}^\infty\beta_k\|y_0\|_k<\frac{\bar r\beta_{k_0+d_1+d_2}}{{a}_{k_0+d_1} }. \end{equation}
\end{center}

Set
\begin{equation}\label{Para2}
 \kappa:=\frac{\sum_{k=0}^\infty\beta_k\|y_0\|_k}{\bar r\alpha_{k_0}}=\frac{f(0,0)}{\bar r\alpha_{k_0}}<\frac{\eta(1-C\gamma-\eta)}{C\gamma+\eta+\eta\varepsilon}.
 \end{equation}
Applying  {Ekeland's variational principle} to the function $f$ on $X\times Y$ endowed with  the distance $d_\varepsilon,$  we may  find $(x_0,z_0)\in X\times Y$ such that
\begin{equation}\label{AE1}
f(x_0,z_0)\le f(0,0),
\end{equation}
\begin{equation}\label{AE2}
d_\varepsilon((x_0,z_0),(0,0))\le \bar r\alpha_{k_0},
\end{equation}
and
\begin{equation}\label{AE3}
f(x,y)+\kappa(d(x,x_0)+\varepsilon d(y,z_0))\ge f(x_0,z_0)\;\;\forall (x,y)\in X\times Y.
\end{equation}
Obviously, $(x_0,z_0)\in\gph F.$ By  relations (\ref{Para1}) and  (\ref{AE2}), 
$$\alpha_{k_0}\min\{r, \|x_0\|_{k_0}\}\le d_\varepsilon((x_0,z_0),(0,0))\le \bar r\alpha_{k_0}<r\alpha_{k_0}. $$
Consequently, $\|x_0\|_{k_0}<r.$ Since according to  (\ref{Cond-Image}), 
$$\beta_{k_0+d_1+d_2}\|y_0\|_{k_0+d_1+d_2}\le f(0,0)<\frac{r\beta_{k_0+d_1+d_2}}{{a}_{k_0+d_1}},$$
it follows $\|y_0\|_{k_0+d_1+d_2}<r/{a}_{k_0+d_1}.$ Furthermore, as
$$\beta_{k_0+d_1+d_2}\|z_0-y_0\|_{k_0+d_1+d_2}\le f(x_0,z_0)\le f(0,0)<\frac{r\beta_{k_0+d_1+d_2}}{{a}_{k_0+d_1}},$$
one obtains
$$\|z_0\|_{k_0+d_1+d_2}\le \|y_0\|_{k_0+d_1+d_2}+ \|z_0-y_0\|_{k_0+d_1+d_2}\le 2r/{a}_{k_0+d_1}.$$ 
From (\ref{AE3}), for all $(x,y)\in\gph F$ one has  
\begin{gather}\label{Iqn-min}
\sum_{k=0}^\infty\beta_k(\|z_0-y_0\|_k-\|y-y_0\|_k)\nonumber\\
\le \eta\sum_{k=0}^\infty\alpha_k(\| x\|_k-\|x_0\|_k)\nonumber\\
+\kappa\left(\sum_{k=0}^\infty\alpha_k\min\{r,\|x-x_0\|_k\}  
+\varepsilon\sum_{k=0}^\infty\beta_k\min\{r,\|y-z_0\|_k\}\right).
\end{gather}
It suffices to prove that $z_0=y_0.$  Assume to  the  contrary that $z_0\not=y_0.$  Setting  $v=y_0-z_0,$ and using the assumption that  $Y$ is standard, there exists a sequence $(v_n)_{n\in \N}$ converging to $v$ such that
\begin{equation}\label{Con-Stand}
\|v_n\|_k\le c_0(v)\|v\|_k,\quad \|v_n\|\le c_1(v_n)^k\;\;\forall n,k\in\N.
\end{equation}
From  condition (ii),  for every $n$, there exists $u_n\in DF(x_0,z_0)^{-1}(v_n)$ such that
\begin{equation}\label{Inver-Contingent}
\|u_n\|_k\le C\left(\frac{\| x_0\|_{k-d_1-d_2}}{m_k{a}_{k-d_1-d_2}}+\nu^\prime_k\right)\|v_n\|_{d_1+d_2}+ {a}_k\|v_n\|_{k+d_2},\;\; \forall k\in\N.
\end{equation}
From  condition (i), for every $n\in\N,$ there exist a real ${c_2(u_n,v_n)}>0$, sequences $t_{n,j}\downarrow 0,$ $u_{n,j}\to u_n,$ $v_{n,j}\to v_n$ as $j\to\infty$ such that
$$(x_0+t_{n,j}u_{n,j},z_0+t_{n,j}v_{n,j})\in\gph F,\;\forall j$$
\begin{equation}\label{Iqn1}
\|v_{n,j}\|_k\le {c_2(u_n,v_n)}(m_k\|u_n\|_{k+d_1}+\|x_0\|_k/{a}_{k-d_2}+\|z_0\|_k+\nu_k),\;\forall j, k,
\end{equation}
\begin{equation}\label{Iqn2}
\|u_{n,j}\|_k\le {c_2(u_n,v_n)}(m_k\|u_n\|_{k+d_1}+\|x_0\|_k+\nu_k),\;\forall j, k.
\end{equation}
Plugging  $x:=x_0+t_{n,j}u_{n,j}$ and $ y:= z_0+t_{n,j}v_{n,j}$ into relation  (\ref{Iqn-min}), one obtains
\begin{gather}\label{Iqn-min-bis}
\sum_{k=0}^\infty\beta_k(\|v\|_k-\|v-t_{n,j}v_{n,j}\|_k)\le \eta\sum_{k=0}^\infty\alpha_k(\| x_0+t_{n,j}u_{n,j}\|_{k+d_1}-\|x_0\|_{k+d_1})\nonumber\\
+\kappa\left(\sum_{k=0}^\infty\alpha_k\min\{r,t_{n,j}\|u_{n,j}\|_k\} +\varepsilon\sum_{k=0}^\infty\beta_k\min\{r,t_{n,j}\|v_{n,j}\|_k\}\right)\;\;\forall n,j.
\end{gather}
We can assume $t_{n,j}\in (0,1),$ for all $n,j.$  Then, one has
$$
\|v-t_{n,j}v_{n,j}\|_k=\|t_{n,j}(v-v_{n,j})+(1-t_{n,j})v\|_k\le t_{n,j}\|v-v_{n,j}\|_k+(1-t_{n,j})\|v\|_k.
$$
It follows that
\begin{equation}
\frac{\|v\|_k-\|v-t_{n,j}v_{n,j}\|_k}{t_{n,j}}\ge \|v\|_k-\|v-v_{n,j}\|_k\;\;\forall n,j,k.\notag
\end{equation}
Combining  this inequality with (\ref{Iqn-min-bis}),  one derives
\begin{equation}\label{@}
\begin{array}{ll}
&\sum_{k=0}^\infty\beta_k(\|v\|_k-\|v-v_{n,j}\|_k)\le \eta\sum_{k=0}^\infty\alpha_k\| u_{n,j}\|_{k+d_1}+\\
&+\kappa\left(\sum_{k=0}^\infty\alpha_k\min\left\{\frac{r}{t_{n,j}},\|u_{n,j}\|_k\right\} +\varepsilon\sum_{k=0}^\infty\beta_k\min\left\{\frac{r}{t_{n,j}},\|v_{n,j}\|_k\right\}\right)\;\;\forall n,j.
\end{array}
\end{equation}
\vskip 0.2cm
\noindent\textit{\textbf{Claim 2.}  For every $n\in\N,$ one has
$$\lim_j\sum_{k=0}^\infty\beta_k(\|v\|_k-\|v-v_{n,j}\|_k)=\sum_{k=0}^\infty\beta_k(\|v\|_k-\|v-v_n\|_k)$$ and
$$\lim_{j\to\infty}\sum_{k=0}^\infty\beta_k\min\left\{\frac{r}{t_{n,j}},\|v_{n,j}\|_k\right\}=\sum_{k=0}^\infty\beta_k\|v_n\|.$$}
\vskip 0.2cm
By relations (\ref{Con-Stand}), (\ref{Inver-Contingent})  and (\ref{Iqn1}), one has
$$\begin{array}{ll}
\|v_{n,j}\|_k\le {c_2(u_n,v_n)}m_k\|u_n\|_{k+d_1}+{c_2(u_n,v_n)}\|x_0\|_k/{a}_{k-d_2}+ {c_2(u_n,v_n)}\|z_0\|_k +{c_2(u_n,v_n)}\nu_k\\
\le {c_2(u_n,v_n)}\left[C\left(\frac{\|x_0\|_{k-d_2}}{{a}_{k-d_2}}+m_k\nu_{k+d_1}^\prime\right)\|v_n\|_{d_1+d_2}+m_k{a}_{k+d_1}\|v_n\|_{k+d_1+d_2}\right]\\+{c_2(u_n,v_n)}\|x_0\|_k/{a}_{k-d_2}+{c_2(u_n,v_n)}\|z_0\|_k\\
\le {c_2(u_n,v_n)}\left[\frac{C}{{a}_{k-d_2}}\|x_0\|_{k-d_2}\|v_n\|_{d_1+d_2}+m_k{a}_{k+d_1}c_1(v_n)^{k+d_1+d_2}\right]\\
+C{c_2(u_n,v_n)}\|v_n\|_{d_1+d_2}m_k\nu_{k+d_1}^\prime +{c_2(u_n,v_n)}\|x_0\|_k/{a}_{k-d_2}+ {c_2(u_n,v_n)}\|z_0\|_k +{c_2(u_n,v_n)}\nu_k.
\end{array}$$
Thus, for every $n, j\in\N,$ 
\begin{equation}\label{conver1}
\begin{array}{ll}
\sum_{k=0}^\infty\beta_k\|v_{n,j}\|_k\\
\le {c_2(u_n,v_n)}\left[C\|v_n\|_{d_1+d_2}\sum_{k=0}^\infty\alpha_{k-d_1-d_2}\|x_0\|_{k-d_2} 
+\sum_{k=0}^\infty \beta_km_k{a}_{k+d_1}c_1(v_n)^{k+d_1+d_2}\right]\\+\sum_{k=0}^\infty\alpha_{k-d_1}\|x_0\|_k+\sum_{k=0}^\infty\beta_k\|z_0\|_k \\
+{c_2(u_n,v_n)}(C\|v_n\|_{d_1+d_2}\sum_{k=0}^\infty\beta_k\nu_{k+d_1}^\prime+\sum_{k=0}^\infty\beta_k\nu_k).
\end{array}
\end{equation}
By (\ref{nu-cond}),
$$C\|v_n\|_{d_1+d_2}\sum_{k=0}^\infty\beta_km_k\nu_{k+d_1}^\prime+\sum_{k=0}^\infty\beta_k\nu_k<+\infty;$$
\tcb{
$$\sum_{k=0}^\infty\beta_km_k{a}_{k+d_1}c_1(v_n)^{k+d_1+d_2}<\infty,$$} 
and in view of relation (\ref{AE1}),
$$\eta\sum_{k=0}^\infty\alpha_k\|x_0\|_{k+d_1}\le f(x_0,z_0)\le f(0,0)<\infty;$$
$$\sum_{k=0}^\infty\beta_k\|z_0\|_k\le \sum_{k=0}^\infty\beta_k\|y_0\|_k+f(0,0)<\infty.$$
Therefore, according to Lebesgue's Dominated Convergence Theorem, relation (\ref{conver1}) yields for every $n\in \N,$
$$\begin{array}{ll}
\displaystyle \lim_{j\to\infty}\sum_{k=0}^\infty\beta_k(\|v\|_k-\|v-v_{n,j}\|_k)
&=\sum_{k=0}^\infty\beta_k(\|v\|_k-\displaystyle\lim_{j\to\infty}\|v-v_{n,j}\|_k)\\
&=\sum_{k=0}^\infty\beta_k(\|v\|_k-\|v-v_{n}\|_k),
\end{array}$$
and
$$\lim_{j\to\infty}\sum_{k=0}^\infty\beta_k\min\left\{\frac{r}{t_{n,j}},\|v_{n,j}\|_k\right\}=\sum_{k=0}^\infty\beta_k\|v_n\|.\quad\quad\qed$$
\vskip 0.1cm
\textit{\textbf{Claim 3.} For every $n\in\N,$ one has
\begin{equation}\label{Conver2}\lim_{j\to\infty}\sum_{k=0}^\infty\alpha_k\| u_{n,j}\|_{k+d_1}=\sum_{k=0}^\infty\alpha_k\| u_{n}\|_{k+d_1},
\end{equation}
and
\begin{equation}\label{Conver3}
\lim_{j\to\infty}\sum_{k=0}^\infty\alpha_k\min\left\{\frac{r}{t_{n,j}},\|u_{n,j}\|_k\right\}=\sum_{k=0}^\infty\alpha_k\| u_{n}\|_{k}.
\end{equation}}
From relations (\ref{Con-Stand}), (\ref{Inver-Contingent}) and (\ref{Iqn2}), for every $k,n,j\in\N,$ one has
$$\begin{array}{ll}\|u_{n,j}\|_k\le {c_2(u_n,v_n)}(m_k\|u_n\|_{k+d_1}+\|x_0\|_k+\nu_k)\\
\le {c_2(u_n,v_n)}\left(C\left(\frac{\|x_0\|_{k-d_2}}{{a}_{k-d_2}}+m_k\nu_{k+d_1}^\prime\right)\|v_n\|_{d_1+d_2}+m_k{a}_{k}\|v_n\|_{k+d_2}+\|x_0\|_k+\nu_k\right)\\
\le {c_2(u_n,v_n)}\left(C\|x_0\|_{k-d_2}\|v_n\|_{d_1+d_2}+m_k{a}_{k}c_1(v_n)^{k+d_2}+\|x_0\|_k\right)\\+ {c_2(u_n,v_n)}\nu_k+{c_2(u_n,v_n)}C\|v_n\|_{d_1+d_2}m_k\nu^\prime_{k+d_1}.
\end{array}.$$
As $\sum_{k=0}^\infty\beta_km_k\nu_{k+d_1}^\prime,$ $\sum_{k=0}^\infty\alpha_k\|x_0\|_{k+d_1}$ and  $\sum_{k=0}^\infty\beta_km_{k+d_1}{a}_{k+d_1}c_1(v_n)^{k+d_1+d_2}$ are convergent series, we deduce (\ref{Conver2}) and (\ref{Conver3}) by Lebesgue's Dominated Convergence Theorem.\quad \quad\qed
\vskip 0.2cm
By virtue of Claims 2 and 3,  by letting $j\to\infty$ in inequality (\ref{@}), one obtains
\begin{equation}\label{@@}
\begin{array}{ll}
&\sum_{k=0}^\infty\beta_k(\|v\|_k-\|v-v_{n}\|_k)\le \eta\sum_{k=0}^\infty\alpha_k\| u_{n}\|_{k+d_1}+\\
&+\kappa\left(\sum_{k=0}^\infty\alpha_k\|u_{n}\|_k +\varepsilon\sum_{k=0}^\infty\beta_k\|v_{n}\|_k\right)\;\;\forall n.
\end{array}
\end{equation} 
Next, using  the first relation of (\ref{Con-Stand}), and the inequalities
$$\sum_{k=0}^\infty\beta_k\|v\|_k\le f(x_0,z_0)\le f(0,0),$$ 
for every $n\in\N,$ one has
$$\sum_{k=0}^\infty\beta_k\|v_n\|_k\le c_0(v)\sum_{k=0}^\infty\beta_k\|v\|_k<\infty.$$
Applying again Lebesgue's Dominated Convergence Theorem, one obtains
\begin{equation}\label{limit1}
\lim_{n\to\infty}\sum_{k=0}^\infty\beta_k(\|v\|_k-\|v-v_{n}\|_k)=\sum_{k=0}^\infty\beta_k(\|v\|_k-\lim_{n\to\infty}\|v-v_{n}\|_k)=\sum_{k=0}^\infty\beta_k\|v\|_k,
\end{equation} 
and
\begin{equation}\label{limit2}
\lim_{n\to\infty}\sum_{k=0}^\infty\beta_k\|v_{n}\|_k=\sum_{k=0}^\infty\beta_k\lim_{n\to\infty}\|v_{n}\|_k=\sum_{k=0}^\infty\beta_k\|v\|_k.
\end{equation}
From (\ref{Inver-Contingent}), one has
$$
\begin{array}{ll}
\sum_{k=0}^\infty\alpha_k\|u_n\|_k\le \sum_{k=0}^\infty\alpha_k\|u_n\|_{k+d_1}\\
\le \sum_{k=0}^\infty C\alpha_k\left(\frac{\| x_0\|_{k-d_2}}{m_{k+d_1}{a}_{k-d_2}}+\nu^\prime_{k+d_1}\right)\|v_n\|_{d_1+d_2}+\sum_{k=0}^\infty \alpha_k{a}_{k+d_1}\|v_n\|_{k+d_1+d_2}\\
\le C \sum_{k=0}^\infty\alpha_k\left(\| x_0\|_{k+d_1}+\nu_{k+d_1}^\prime\right)\|v_n\|_{d_1+d_2}+\sum_{k=0}^\infty \beta_{k+d_1+d_2}\|v_n\|_{k+d_1+d_2}.
\end{array}
$$
As
$$\eta\sum_{k=0}^\infty\alpha_k\| x_0\|_{k+d_1}\le f(x_0,z_0)\le f(0,0)=\sum_{k=0}^\infty\beta_k\|y_0\|_k,$$
one deduces that
\begin{equation}\label{Estim1}
\begin{array}{ll}
\sum_{k=0}^\infty\alpha_k\|u_n\|_k&\le \sum_{k=0}^\infty\alpha_k\|u_n\|_{k+d_1}\\
&\le \frac{C\gamma}{\eta} \sum_{k=d_1+d_2}^\infty\beta_k\|v_n\|_{d_1+d_2}+\sum_{k=0}^\infty \beta_{k+d_1+d_2}\|v_n\|_{k+d_1+d_2}\\
&\le \frac{C\gamma}{\eta} \sum_{k=d_1+d_2}^\infty\beta_k\|v_n\|_{k}+\sum_{k=0}^\infty \beta_{k+d_1+d_2}\|v_n\|_{k+d_1+d_2}\\
&\le \left(\frac{C\gamma}{\eta}+1\right)\sum_{k=0}^\infty\beta_k\|v_n\|_k.
\end{array}
\end{equation}
By virtue of this inequality, letting $n\to\infty$ in relation (\ref{@@}), in view of relations (\ref{limit1}) and (\ref{limit2}), one obtains
$$\sum_{k=0}^\infty\beta_k\|v\|_k\le \left[(C\gamma+\eta)+\kappa(C\gamma/\eta+1+\varepsilon)\right]\sum_{k=0}^\infty\beta_k\|v\|_k,$$
from which  {it} follows that
$$\kappa\ge \frac{\eta(1-C\gamma-\eta)}{C\gamma+\eta+\eta\varepsilon}.$$
This contradicts (\ref{Para1}), (\ref{Para2}). The proof is completed.\hfill{$\Box$}
\vskip 0.5cm
\begin{remark} 
We can add to the conclusion of  Theorem  \ref{MS},  that the inverse image  $x$ of   an arbitrary $y\in Y$ is controlled by the distance to the reference point $\bar x$. More precisely thanks to   relation (\ref{AE2}) we have $d(x,\bar x)< r\alpha_{k_0},$ where $d$ is a metric defining $X.$ 
\end{remark}   
\vskip 3mm
In the last part of this section, we consider two graded Fr\'echet spaces $X$ and  $Y$, a topological space $P$ and  a multifunction $F:X\times P\rightrightarrows Y$. For $p\in P,$ set
\begin{equation}\label{Sol-Map}
S(p)=\{x\in X:\;\; 0\in F(x,p)\}.
\end{equation}
The multifunction $S:P\rightrightarrows X$ is {referred to as}  the {\it solution mapping} associated to $F.$ For $p\in P,$ denote by $F_p:=F(\cdot,p):X\rightrightarrows Y.$ {By making use }Theorem \ref{MS}, we derive the following implicit multifunction theorem.  
\begin{theorem}\label{IFT}
Let $X$, $Y$ be graded Fr\'{e}chet spaces and let $P$ be a topological space. Consider a multifunction $F:X\times P\rightrightarrows Y$ and {a given point}  $(\bar x,\bar p)\in X\times P$ with $0\in F(\bar x,\bar p)$. Assume furthermore that $Y$ is standard. Suppose also that there are integers $k_0, d_1,d_2,$ real numbers $r_0\in(0,+\infty]$, $C>0$  and non-decreasing sequences of non-negative   reals $(\mu_k)_{k\in\N},$ $(\mu^\prime_k)_{k\in\N},$ and $(m_k)_{k\in \N},$  $({a}_k)_{k\in \N}$ with $m_k\ge 1,$ ${a}_k\ge 1$ such that  the following conditions are satisfied:
\begin{itemize}
\item[(i)] For each $p\in P,$ $\gph F_p$ is closed and the multifunction $F(\bar x,\cdot):P\rightrightarrows Y$ is {lower semicontinuous } at $\bar p;$ 
\item[(ii)]  For all $p$ near $\bar p,$ for all $(x,y)\in \gph F_p$  with $x\in B_X(\bar x,k_0,r_0),$ $y\in B_Y(0,k_0+d_1+d_2,2r_0/{a}_{k_0+d_1}),$ for every  $v\in DF_p(x,z)u$ with $u\in X,$  there exist {$c_2(u,v)>0$} and  sequences $t_n\downarrow 0,$ $ u_n\to u$ and $v_n\to v$ with \\$(x+t_nu_n,y+t_nv_n)\in\gph F_p$, such that for all $n\in\N,$ all $k\in\N,$
$$\|v_n\|_k\le {c_2(u,v)}(m_k\|u\|_{k+d_1}+\|x-\bar x\|_k/{a}_{k-d_2}+\|y\|_k+\mu_k) $$
and 
$$\|u_n\|_k\le {c_2(u,v)}(m_k\|u\|_{k+d_1}+\|x-\bar x\|_k+\mu_k) ;$$
\item[(ii)]  For $p\in P$ near $\bar p,$ for all $(x,y)\in\gph F_p$ with  $x\in B_X(\bar x,k_0,r_0),$ $y\in B_Y(0,k_0+d_1+d_2,2r/{a}_{k_0+d_1}),$  for  every $v\in Y,$  there exists $u\in DF_p(x,y)^{-1}(v)$ such that
$$\|u\|_k\le \frac{C}{m_k{a}_{k-d_1-d_2}}\|x-\bar x\|_{k-d_1-d_2}\|v\|_{d_1+d_2}+{a}_k\|v\|_{k+d_2},\;\; \forall k\in\N.$$
\end{itemize}
Then for every $\tau>{a}_{k_0+d_1}$, there exist $r\in (0,r_0)$ and a neighborhood $W$ in $P$ of $\bar p$ such that
\begin{equation}\label{MIq}
d_{k_0}(x,S(p))\le \tau d_{k_0+d_1+d_2}(0, F(x,p)),
\end{equation}
{for all $(x,p)\in B_X(\bar x, k_0,r)\times W.$} 
\end{theorem}
\vskip 0.2cm
\textit{Proof.} Pick a positive real $r\in (0,r_0)$ such as $Cr<1/2.$ Since $F(\bar x,\cdot)$ is lower semicontinuous, for any $\varepsilon>0$, there exists a neighborhood $W$ of $\bar p$ in $P$ such that $F(\bar x,p)\cap B_Y(0, k_0+d_1+d_2,\varepsilon)\not=\emptyset,$ for all $p\in W.$ That is, 
$$d_{k_0+d_1+d_2}(0, F(\bar x,p))<\varepsilon,\quad\forall p\in W.$$
 Suppose that for this neighborhood $W,$  conditions (i) and (ii) are satisfied for all $p\in W.$ For a given  $p\in W,$ there is $\bar y_p\in F(\bar x,p)$ such as
 $$\|\bar y_p\|_{k_0+d_1+d_2}\le (1+\varepsilon)d_{k_0+d_1+d_2}(0, F(\bar x,p))<\varepsilon.$$ 
 Let $\tau>{a}_{k_0+d_1}$ and $p\in W$ with $0\notin F(\bar x,p)$ be given. Pick $\varepsilon>0$ with $\tau\varepsilon<r/4$ and ${a}_{k_0+d_1}<\tau/(1+\varepsilon)$ and $0<r^\prime<r$ with 
 $${a}_{k_0+d_1}\|\bar y_p\|<r^\prime<\frac{\tau}{1+\varepsilon}\|\bar y_p\|_{k_0+d_1+d_2}\le \tau
 d_{k_0+d_1+d_2}(0, F(\bar x,p)).$$
 Take a sequence $(\beta_k)_{k\in \N}$ of non-negative reals with unbounded support such that
{\begin{equation}\label{beta-seq@bis}
\begin{array}{ll}\sum_{k=0}^\infty\beta_k(\|\bar y_p\|_k+\mu_k)<+\infty
\quad \mbox{and}\quad \sum_{k=0}^\infty\beta_km_k{a}_{k+d_1}n^k<+\infty;
\end{array}
\end{equation}}
\begin{equation}\label{Cond-image@}
\frac{C\sum_{k=0}^\infty\beta_k\|\bar y_p\|_k}{\sum_{k=d_1+d_2}^\infty\beta_k}<1,
\end{equation}
and
\begin{equation}\label{Cond-Image@}
\sum_{k=0}^\infty\beta_k\|y_p\|_k\left(1-\sqrt{\frac{C\sum_{k=0}^\infty\beta_k\|y_p\|_k}{\sum_{k=d_2}^\infty\beta_k}}\right)^{-2}<\frac{r^{\prime}\beta_{k_0+d_1+d_2}}{{a}_{k_0+d_1}}.
\end{equation}
It is not difficult to show the existence of such a sequence $(\beta_k)_{k\in \N}$. Let us apply Theorem \ref{MS} with $F(\cdot,p),$ $r^\prime$, $\bar y_p$ and $0,$ instead of $F,$ $r,$ $\bar y$ and $y,$ respectively. Obviously, $(i)$ and (ii) are satisfied with $\nu_k=\mu_k+\|\bar y_p\|_k$, $\nu_k^\prime=0.$ We can find $\bar x_p\in S(p)$ such that
\begin{align}\label{sol1}
d_k(\bar x,S(p))&\le\|\bar x_p-\bar x\|_{k_0}\nonumber\\
&\le r^\prime 
<\tau
 d_{k_0+d_1+d_2}(0, F(\bar x,p)).&
\end{align} 
Now let  $x\in B_X(\bar x,k_0,r/2)$ and $p\in W$. If $d_{k_0+d_1+d_2}(0,F(x,p))\ge r/(2\tau)$ {be  given}, then by relation (\ref{sol1}), there is $\bar x_p\in S(p)$ such that $\|\bar x-\bar x_p\|_{k_0}\le \tau\varepsilon<r/4.$ Therefore, 
\begin{align}\label{sol2}
d_{k_0}(x,S(p))\le\|x-x_p\|_{k_0}&\nonumber\\
&\le \|x-\bar x\|_{k_0}+\|\bar x-x_p\|_{k_0}\nonumber\\
&\le r/2+r/4<\tau d_{k_0+d_1+d_2}(0,F(x,p)).
\end{align}
Let us consider the case $d_{k_0+d_1+d_2}(0,F(x,p))<r/(2\tau).$ Then there is $y_p\in F(x,p))$ with $\|y_p\|_{k_0+d_1+d_2}<r/2.$ Pick a non-negative  number $r^\prime<r/2$ with
$${a}_{k_0+d_1}\|y_p\|_{k_0+d_1+d_2}<r^\prime<\tau\|y_p\|_{k_0+d_1+d_{2}.}$$
Since $C(\|y_p\|_{k_0+d_1+d_2}+\|x-\bar x\|_{k_0})<Cr<1/2,$ we can pick a sequence $(\beta_k)_{k\in \N}$ of non-negative reals with unbounded support such that
$$\sum_{k=0}^\infty\beta_k(\|y_p\|_k+\mu_k)<+\infty\quad\mbox{and}\quad\sum_{k=0}^\infty\beta_km_k\|x-\bar x\|_{k-d_2}<+\infty;$$
$$\sum_{k=0}^\infty\beta_k<+\infty\quad\mbox{and}\quad \sum_{k=0}^\infty\beta_km_k{a}_{k+d_1}n^k<+\infty\quad\mbox{and}\quad  {(\ref{nu-cond})};$$
$$
\frac{C\sum_{k=0}^\infty\left(\beta_k\|y_p\|_k+\beta_{k+d_1+d_2}\|x-\bar x\|_{k-d_2}\right)}{\sum_{k=d_1+d_2}^\infty\beta_k}:=s<1;
$$
$$
\sum_{k=0}^\infty\beta_k\|y_p\|_k\left(1-\sqrt{s}\right)^{-2}<\frac{r^\prime\beta_{k_0+d_1+d_2}}{{a}_{k_0+d_1}}.
$$
Then apply Theorem \ref{MS} with $F(\cdot,p),$ $r^\prime$, $y_p$ $0,$ $x$ instead of $F,$ $r,$ $\bar y$ $y,$ and $\bar x,$ respectively, since conditions $(i)$ and $(ii)$ are verified for $\nu_k=\|x-\bar x\|_k+\|y_p\|_k+\mu_k$, $\nu_k^\prime=\|x-\bar x\|_{k-d_1-d_2}.$ We  obtain the existence of $x_p\in S(p)$ verifying 
$$
d_k(x,S(p))\le\|x_p-x\|_{k_0}\le r^\prime<\tau
 d_{k_0+d_1+d_2}(0, F(\bar x,p)).
$$
Thus (\ref{MIq}) is shown.\hfill{$\Box$}
\vskip 0.5cm
To a given   multifunction $F:X\rightrightarrows Y,$  we associate the new multifunction $\varPhi: X\times Y\rightrightarrows Y$ defined by 
$$\varPhi(x,y)=F(x)-y,\quad (x,y)\in X\times Y.$$Applying Theorem \ref{IFT} to the multifunction  $\varPhi$,
we derive the following result of metric regularity type  in Fr\'{e}chet spaces.
\begin{corollary}\label{MR}
Let $X,Y$ be graded Fr\'{e}chet spaces and assume that $Y$ is standard. Let $F:X\rightrightarrows Y$ be a closed multifunction and let  $(\bar x,\bar y)\in\gph F$ {be  given}. Suppose also that there are integers $k_0, d_1,d_2,$ real numbers $r\in(0,+\infty]$, $C>0$  and non-decreasing sequences of non-negative reals $(\nu_k)_{k\in \N},$ $(\nu^\prime_k)_{k\in \N},$ and $(m_k)_{k\in \N},$  $({a}_k)_{k\in \N}$ with $m_k\ge 1,$ ${a}_k\ge 1$ such that  the { conditions (i) and (ii)} of Theorem \ref{MS} are satisfied. Then for every $\tau>{a}_{k_0+d_1}$, there exists a neighborhood $W$ in $Y$ of $\bar y$ such that
\begin{equation}\label{MRI}
d_{k_0}(x,F^{-1}(y))\le \tau d_{k_0+d_1+d_2}(y, F(x)),
\end{equation}
for all $(x,y)\in B_X(\bar x, k_0,r)\times W.$
\end{corollary}
\vskip 0.2cm

When the mapping $F: X\rightarrow Y$ is G\^{a}teaux differentiable, Theorem \ref{MS} yields the following theorem which covers Theorem 3 in \cite{E-IHP}.
\begin{theorem}\label{MS-Diff}
Let $F:X\rightarrow Y$ be a continuous mapping between graded Fr\'{e}chet spaces $X,Y$ with $Y$ standard, and let $\bar x\in X$ be given. Suppose $F$ is G\^{a}teaux differentiable on $X$ with derivative $DF(x)$, and that there are integers $k_0, d_1,d_2,$ real numbers $r\in(0,+\infty]$, $C>0$  and non-decreasing sequences of non-negative reals $(\nu_k)_{k\in \N},$ and $(m_k)_{k\in \N},$  $({a}_k)_{k\in \N}$ with $m_k\ge 1,$ ${a}_k\ge 1$ such that  the following conditions are satisfied:
\begin{itemize}

\item[{(i)}] For every $u\in X$, there exists $c_2(u)>0$ such that for all $x\in B_X(\bar x,k_0,r),$ all $k\in\N,$
$$\|DF(x)u\|_k\le c_2(u)(m_k\|u\|_{k+d_1}+\|x-\bar x\|_k/{a}_{k-d_2}+\|F(x)-F(\bar x)\|_k+\nu_k) $$
\item[{(ii)}]  For all $x\in B_X(\bar x,k_0,r),$ there exists a linear mapping $L(x): Y\rightarrow X$ such that $DF(x)L(x)=I_Y$ (the identity mapping on $Y$) and for every $v\in Y,$  
$$\|L(x)v\|_k\le \frac{C}{m_k{a}_{k-d_1-d_2}}\|x-\bar x\|_{k-d_1-d_2}\|v\|_{d_1+d_2}+{a}_k\|v\|_{k+d_2},\;\; \forall k\in\N.$$
\end{itemize}
Let $(\beta_k)_{k\in \N}$ be a sequence of non-negative reals with unbounded support such that
\begin{equation}
\sum_{k=0}^\infty\beta_k\nu_k<+\infty,
\end{equation}
Then, for every $y\in Y$ with
\begin{equation}
\frac{C\sum_{k=0}^\infty\beta_k\|y-F(\bar x)\|_k}{\sum_{k=d_1+d_2}^\infty\beta_k}<1,
\end{equation}
and
\begin{equation}
\sum_{k=0}^\infty\beta_k\|y-F(\bar x)\|_k\left(1-\sqrt{\frac{C\sum_{k=0}^\infty\beta_k\|y-F(\bar x)\|_k}{\sum_{k=d_1+d_2}^\infty\beta_k}}\right)^{-2}<\frac{r\beta_{k_0+d_1+d_2}}{{a}_{k_0+d_1}},
\end{equation}
there exists $x\in B_X(\bar x,k_0,r)$ such that $F(x)=y.$
\end{theorem}
\vskip 0.2cm
\textit{Proof.} It suffices to show that conditions (i) and (ii) of Theorem \ref{MS} are satisfied. Indeed, for every $x\in B(\bar x, k_0,r),$ $v\in Y,$ setting $u=L(x)v,$ one has obviously $u\in DF^{-1}(x)(v).$ Thus condition (2) implies (ii) of Theorem \ref{MS}. For given $u\in X$ and $x\in B(\bar x,k_0,r),$ pick $\bar t\in (0,1)$ such that $\|x+tu\|_{k_0}<r$ for all $t\in [0,\bar t].$ For each $k\in\N$, define the function $f_k:[0,\bar t]\to\R$ by
$$f_k(t)=\|F(x+tu)-F(\bar x)\|_k,\;\; t\in [0,t_0].$$
Obviously, $f$ has a right derivative everywhere and 
$$f_{k+}^\prime(t)\le \|DF(x+tu)\|_k,\quad\forall t\in[0,\bar t).$$
Therefore, by assumption {$(i)$}, one has
$$f^\prime_{k+}(t)\le c_2(u)(m_k\|u\|_{k+d_1}+\|x+tu-\bar x\|_k/{a}_{k-d_2}+f(t)+\nu_k),\;\;\forall t\in[0,\bar t),$$
and consequently,
$$f^\prime_{k+}(t)-c_2(u)f_k(t)\le c_2(u)((m_k+1)\|u\|_{k+d_1}+\|x-\bar x\|_k+\nu_k),\;\;\forall t\in[0,\bar t).$$
Equivalently,
\[
\begin{array}{ll}
e^{-tc_2(u)}[f^\prime_{k+}(t)-c_2(u)f_k(t)]&\\
\le e^{-tc_2(u)}c_2(u)((m_k+1)\|u\|_{k+d_1}+\|x-\bar x\|_k/{a}_{k-d_2}+\nu_k),\;\;\forall t\in[0,\bar t). &$$
\end{array}
\]
By integration, one obtains
\[
\begin{array}{ll}
e^{-tc_2(u)}f_k(t)-f_k(0)&\\
\le (1-e^{-tc_2(u)})((m_k+1)\|u\|_{k+d_1}+\|x-\bar x\|_k/{a}_{k-d_2}+\nu_k).&
\end{array}
\]
That is, for all $t\in [0,\bar t),$
\[
\begin{array}{ll}
m_k\|u\|_{k+d_1}+\|F(x+tu)-F(\bar x)\|_k\ &\\
\le e^{tc_2(u)}[(m_k+1)\|u\|_{k+d_1}
+ \|x-\bar x\|_k/{a}_{k-d_2}+\|F(x)-F(\bar x)\|_k+\nu_k].&\\
\end{array}
\]
This together  with $({i})$  yield
\[
\begin{array}{ll}
\|DF(x+tu)u\|_k&\\
\le c_2(u)e^{tc_2(u)}[(m_k+1)\|u\|_{k+d_1}+\|x-\bar x\|_k/{a}_{k-d_2}+\|F(x)-F(\bar x)\|_k+\nu_k],\;\;\forall t\in [0,t_0).&
\end{array}
\]
Next, pick a sequence $(t_n)_{n\in\N}$ converging to $0,$ with $t_n\in (0,\bar t),$ and set
$$u_n:=u,\quad v_n:=\frac{F(x+t_nu)-F(x)}{t_n},\;\; n\in\N.$$
Then,  $F(x)+t_nv_n=F(x+t_nu_n)$ and $\lim_{n\to\infty} (u_n,v_n)=(u,DF(x)u).$ Setting $C_2(u)=c_2(u)e^{{\bar t} c_2(u)}$, the Mean Value inequality yields, 
\[
\begin{array}{ll}
\|v_n\|_k&\\
\le \sup_{s\in [0,t_n]}\|DF(x+su)u\|_k&\\
\le C_2(u)[(m_k+1)\|u\|_{k+d_1}+\|x-\bar x\|_k/{a}_{k-d_2}+\|F(x)-F(\bar x)\|_k+\nu_k], \text{for all} \;n\in \N,&
\end{array}
\]
Thus condition $(ii)$ of Theorem \ref{MS} holds.\hfill{$\Box$}
\vskip 0.5cm
\begin{corollary}\label{MRD}
Under the assumptions of Theorem \ref{MS-Diff}, for every $\tau>{a}_{k_0+d_1}$, there exist $r>0$ and a neighborhood $W$ in $Y$ of $\bar y$ such that
\begin{equation}\label{MRID}
d_{k_0}(x,F^{-1}(y))\le \tau \|y- F(x)\|_{k_0+d_1+d_2},
\end{equation}
for all $(x,y)\in B_X(\bar x, k_0,r)\times W.$
\end{corollary}
\vskip 0.5cm
Using the same argument as in the proof of Theorem \ref{MS-Diff}, from Theorem \ref{IFT}, we obtain, when $F(\cdot,p)$ is G\^{a}teaux differentiable,  the following   implicit multifunction theorem for the system (\ref{Sol-Map}). 
\begin{corollary}\label{IFT-Diff} Let $X$, $Y$ be graded Fr\'{e}chet spaces and let $P$ be a topological space. Consider a mapping $F:X\times P\rightarrow Y$ and {a given point} $(\bar x,\bar p)\in X\times P$ with $F(\bar x,\bar p)=0$. Assume furthermore that $Y$ is standard. Suppose $F(\cdot,p)$ is G\^{a}teaux differentiable on $X$ with derivative $DF_p(x)$ for $p$ near $\bar p,$ and that there are integers $k_0, d_1,d_2,$ real numbers $r_0\in(0,+\infty]$, $C>0$  and non-decreasing sequences of non-negative  reals $(\nu_k)_{k\in \N},$  and $(m_k)_{k\in \N},$  $({a}_k)_{k\in \N}$ with $m_k\ge 1,$ ${a}_k\ge 1$ such that  the following conditions are satisfied:
\begin{itemize}
\item [(i)] The mapping $F$ is continuous at $(\bar x,\bar p);$

\item[(ii)] For $p$ near $\bar p,$ for every $u\in X$, there exists {$c_2(u)>0$} such that for all $x\in B_X(\bar x,k_0,r_0),$ all $k\in\N,$
$$\|DF_p(x)u\|_k\le {c_2(u)}(m_k\|u\|_{k+d_1}+\|x-\bar x\|_k/{a}_{k-d_2}+\|F_p(x)-F_p(\bar x)\|_k+\nu_k); $$

\item[(iii)] For $p$ near $\bar p$ and for all $x\in B_X(\bar x,k_0,r_0),$ there exists a linear mapping $L_p(x): Y\rightarrow X$ such that $DF_p(x)L_p(x)=I_Y$ (the identity mapping on $Y$) and for every $v\in Y,$  
$$\|L_p(x)v\|_k\le \frac{C}{m_k{a}_{k-d_1-d_2}}\|x-\bar x\|_{k-d_1-d_2}\|v\|_{d_1+d_2}+{a}_k\|v\|_{k+d_2},\;\; \forall k\in\N.$$
\end{itemize}
Then,  for every $\tau>{a}_{k_0+d_1}$, there exist $r>0$ and a neighborhood $W$ in $P$ of $\bar p$ such that
\begin{equation}\label{MIq-bis}
d_{k_0}(x,S(p))\le \tau \|F(x,p)\|_{k_0+d_1+d_2},
\end{equation}
for all $(x,p)\in B_X(\bar x, k_0,r)\times W.$
\end{corollary}
\section{Application: Differential equations in Fr\'{e}chet spaces}
 In this final section, we present an application 
to the existence of solutions for ordinary differential equations in Fr\'{e}chet spaces.
Let $X$ be a graded Fr\'{e}chet space, let $U\subseteq X$ be an open set, let $t_0\in\R$, $r_0>0$ be given, and let $f:[t_0-r_0,t_0+r_0]\times U\to X$ be a continuous mapping. For given $r>0$ and $x_0\in U$, consider the initial value problem:
$$(DEF)\quad\quad\quad\left\{\begin{array}{ll}
x^\prime(t)&=f(t,x(t)),\quad t\in [t_0-r,t_0+r],\\
x(t_0)&=x_0.
\end{array}\right.$$
 When the function $f$   is of class  $C^2$, Poppenberg \cite{Pop-Studia}
 established a result on existence of solutions for equation (DEF). In the following theorem, the data function is assumed merely to be G\^{a}teaux differentiable. \vskip 0.2cm
\begin{theorem}\label{ODEF} Let $E$ be graded Fr\'{e}chet space such that $X$ and $C([-1,1],E)$ are standard. Suppose that the function  $f$ is continuous on  $[t_0-r_0,t_0+r_0]\times U.$ Suppose also that  for each $t\in [t_0-r_0,t_0+r_0],$ $f(t,\cdot): U\to X$ is G\^{a}teaux differentiable on $U$ with derivative $D_xf(t,\cdot),$ and that there is a non-decreasing sequence of non-negative  reals $(c_k)_{k\in\N}$ such that for all $(t,x)\in [t_0-r_0,t_0+r_0]\times U,$ one has
\begin{equation}\label{Derivate-Cond}
\|D_xf(t,x)u\|_k\le c_k\|u\|_{k},\quad\mbox{for all}\;u\in X,\; k\in\N.
\end{equation}
Then,  there is $r\in (0,r_0]$ such that problem {(DEF)}  has a solution $x(\cdot)\in C^1([t_0-r,t_0+r],E).$ If   in addition,  $f$ is a $C^1-${mappings}  on $[t_0-r_0,t_0+r_0]\times U$, and that, say, for a sequence  $(c_k)_{k\in\N}$ above,
\begin{equation}\label{df-bounded}
\|Df(t,x)\|_k\le c_k\;\; \forall k\in\N,\;\forall (t,x)\in [t_0-r_0,t_0+r_0]\times U,
\end{equation}
 then the solution $x(\cdot)$ is unique.
\end{theorem}
\vskip 0.2cm
\textit{Proof.} Using the transformations $t=t_0+rs,$ $z(s)=x(t_0+rs)-x_0,$ $s\in[-1,1],$  we can rewrite   problem (DEF) as
$$(DEF1)\quad\quad\quad\left\{\begin{array}{ll}
z^\prime(s)&=rf(t_0+rs,z(s)+x_0),\quad s\in [-1,1],\\
z(0)&=0.
\end{array}\right.$$
Denote by
$$W=\{(z,r)\in C^1([-1,1],X)\times \R: \;\;r\in (-r_0,r_0),\; z(s)\in U\;\forall s\in [-1,1] \}. $$ $W$ is an  open subset of the graded Fr\'{e}chet space $C^1([-1,1],X)\times \R.$ 
Set $F: W\rightarrow C([-1,1],X)\times X,$ defined by
$$F(z,r)=(z^\prime(s)-rf(t_0+rs,z(s)+x_0), z(0)),\;s\in [-1,1],\; (z,r)\in W.$$
Then for each $r\in (0,r_0),$ $F(r,\cdot)$ is G\^{a}teaux differentiable on $W$ with  derivative given by
$$D_zF(z,r)u=(u^\prime(s)-rD_xf(t_0+rs,z(s)+x_0)u(s), u(0)),\; s\in [-1,1],\; (r,s)\in W,\; u\in C^1([-1,1],X).$$
Obviously, $(0,0)\in W$ and $F(0,0)=(0,0),$ and moreover, $z\in C^1([-1,1],X)$ is a solution of problem (DEF1) with respect to $r\in (-r_0,r_0)$ if and only if
$F(z,r)=(0,0).$
So it suffices to show that the mapping $F$ verifies all the assumptions $(i), (ii)$ and $(iii)$ of Corollary \ref{IFT-Diff} with $(\bar x,\bar p):=(0,0),$ $(x,p):=(x,r).$ Assumption $(i)$ is obvious. To verify $(ii)$, for any $k\in \N,$ for $(z,r)\in W,$ $u\in C^1([-1,1],X),$ making use of relation (\ref{Derivate-Cond}), one has
$$\begin{array}{ll}
\|D_zF(z,r)u\|_k&=\|u^\prime(\cdot)-rD_xf(t_0+rs,z(\cdot)+x_0)u(\cdot)\|_k+\|u(0)\|_k\\
&\le \|u^\prime(\cdot)\|_k +r_0c_k\|u(\cdot)\|_k +\|u(\cdot)\|_k\le (2+r_0c_k)\|u\|_k.
\end{array}$$
Thus $(ii)$  follows. To verify $(iii)$, for $(z,r)\in W,$ $(v,v_0)\in C([-1,1],X)\times X,$ then $D_zF(z,r)u=(v,v_0),$ $u\in C^1([-1,1],X)$  if and only if $u$ is a solution of the linear differential equation:
$$(LDE)\quad\quad\quad\left\{\begin{array}{ll}
u^\prime(s)&=A(s)u(s)+v(s),\quad s\in [-1,1],\\
u(0)&=v_0
\end{array}\right.$$
where, $A(s):=rD_xf(t-0+rs, z(s)+x_0)$ ($s\in [-1,1]$) is a continuous linear mapping from $X$ to itself, according to condition (\ref{Derivate-Cond}). Thanks to Proposition 3.4 in \cite{Pop-Studia}, Problem (LDE) has a unique solution $u\in C^1([-1,1],X).$ That is, $D_zF(z,r)$ is invertible with $L(z,r):=D_zF(z,r)^{-1}(v,v_0)=u,$ $u$ solving (LDE), for $(v,v_0)\in C([-1,1],X)\times X.$ Let now $(v,v_0)\in C([-1,1],X)\times X,$ be given, and let $u\in C^1([-1,1],X)$ be a solution of problem (LDE). One has
\begin{equation}\label{A-estim}
\|A(s)u(s)\|_k\le r_0c_k\|u(s)\|_k,\;\;\mbox{for all}\;k\in\N,\; s\in [-1,1]. 
\end{equation}
Therefore, for any $k\in\N,$  by considering the corresponding integral equation of (LDE), one has
$$\|u(t)\|_k\le \|v\|_k+\|v_0\|_k+r_0c_k\int_0^t\|u(s)\|_k,\;\;\forall t\in[0,1].$$ 
Thanks to the Gronwall lemma, applied to the function $\alpha(t):=\|u(t)\|_k$, $t\in[0,1],$ one obtains 
$$\|u(t)\|_k\le e^{r_0c_kt}\|(v,v_0)\|_k,\;\;\mbox{for all}\;\; t\in[0,1].$$
Similarly, for $t\in [-1,0]$, by setting $w(s):=u(-s),$ $s\in [0,1],$ one also has
$$\|w(t)\|_k\le \|v\|_k+\|v_0\|_k+r_0c_k\int_0^t\|w(s)\|_k,\;\;\forall t\in[0,1].$$
Hence,
$$\|u(t)\|_k\le e^{-r_0c_kt}\|(v,v_0)\|_k,\;\;\mbox{for all}\;\; t\in[-1,0].$$
Thus, 
\begin{equation}\label{u-estim}
\sup_{s\in[-1,1]}\|u(s)\|_k\le e^{r_0c_k}\|(v,v_0)\|_k,\quad\forall k\in\N.
\end{equation}
Furthermore, from  equation (LDE),  using relations (\ref{A-estim}), (\ref{u-estim}), one obtains
$$\begin{array}{ll}\sup_{s\in[-1,1]}\|u^\prime(s)\|_k&\le r_0c_k\sup_{s\in[-1,1]}\|u(s)\|_k +\sup_{s\in[-1,1]}\|v(s)\|_k\\
&\le r_0c_ke^{r_0c_k}\|(v,v_0)\|_k +\|v\|_k\le (r_0c_ke^{r_0c_k}+1)\|(v,v_0)\|_k.
\end{array}$$ 
Hence,
$$\|u(\cdot)\|_k\le (e^{r_0c_k}+r_0c_ke^{r_0c_k}+1)\|(v,v_0)\|_k, $$
and $(iii)$ follows. According to Corollary \ref{IFT-Diff}, there is a neighborhood $U(x_0)\subseteq U$ of $x_0$ and $r_1>0$, as well as a sequence of non-negative  reals $(\tau_k)$  such that for all $k\in \N,$
$$d_{k}(x,S(r))\le \tau_k\|F(x,r)\|_{k},\;\;\mbox{for all}\;\; (x,r)\in U(x_0)\times [0,r_1],$$
where, $$S(r)=\{z\in C^1([-1,1],X):\;\; (z,r)\in W,\; F(z,r)=0\}.$$
Suppose now $f$ is a mapping of class $C^1$on $[t_0-r_0,t_0+r_0]\times U.$ We shall show that when $r>0$ is sufficiently small, there is a unique $z\in C^1([-1,1],X)$   with $(z,r)\in W$ such that $F(x,r)=(0,0).$ Indeed,
since $f$ is a $C^1$-mapping on $[t_0-r_0,t_0+r_0]\times U,$ there is a sequence of non-negative  reals $(\varepsilon_k)$ such that for all $k\in\N,$
\[
\begin{array}{ll}
\|f(t_1,x_1)\|_k-f(t_2,x_2)\|_k&\\
\le\|Df(t_0,x_0)\|_k\|(t_1-t_2,x_1-x_2)\|_k&\\
+\varepsilon_k\|(t_1-t_2,x_1-x_2)\|_k,\;\;\forall (t_1,x_1),(t_2,x_2) \in [t_0-r_0,t_0+r_0]\times U.
\end{array}
\]
This shows that there is some sequence of non-negative  $(d_k)$ such that
\begin{equation}\label{f-bounded}
\begin{array}{ll}
\|f(t_1,x_1)\|_k-f(t_2,x_2)\|_k\le& d_k\|(t_1-t_2,x_1-x_2)\|_k\\
&\mbox{for all}\;\; k\in\N,\; (t,x)\in [t_0-r_0,t_0+r_0]\times U.
\end{array}
\end{equation}
Let  $r\in (0,r_0)$ {be given,} and let {$z_{1}, z_{2}$ in $C^1([-1,1],X)$ be}  two solutions of equation (DEF1) with respect to $r.$ Then, one has
$$\left\{\begin{array}{ll}
z_1^\prime(s)-z_2^\prime(s)&=r[f(t_0+rs,z_1(s)+x_0)-f(t_0+rs,z_2(s)+x_0),\quad s\in [-1,1],\\
z_1(0)-z_2(0)&=0.
\end{array}\right.$$
By making use similarly of the Gronwall lemma as above for this differential equation  with its solution $z_1-z_2,$ one derives $z_1=z_2.$ 
The proof is completed.\hfill{$\Box$}
\begin{remark} \rm{ Note that, by the same argument, the conclusion of the preceding theorem remains   true if condition (\ref{Derivate-Cond}) is replaced by the following one:

\noindent there is a constant $C>0$ such that for all $k\in\N,$ one has
\begin{equation}\label{Derivate-Cond-bis}
\|D_xf(t,x)u\|_k\le C(\|x\|_k+\|u\|_{k}),\quad\mbox{for all}\;u\in X,\; (t,x)\in [t_0-r_0,t_0+r_0]\times U.
\end{equation} }
\end{remark}

 \section{Acknowledgements}

The authors are grateful to the referee who has made some constructive comments that have helped to improve the presentation of the paper.

   \bibliographystyle{model1b-num-names}
  \bibliography{Ekeland2.bib}
\end{document}